\documentclass{article} 
\title{Finite subgraphs of uncountably chromatic graphs}

\author{P\'eter Komj\'{a}th\thanks{Research  partially supported by Hungarian 
                                   National Research Grant T 032455.}
                                   \\
         Saharon Shelah\thanks{This research was supported by the 
         Israel  Science Foundation. Publication 788.}} 
                                   
\renewcommand\a{\alpha}
\renewcommand\b{\beta}
\newcommand\g{\gamma}
\renewcommand\d{\delta}

\renewcommand\th{\theta}
\renewcommand\k{\kappa}
\renewcommand\l{\lambda}
\renewcommand\o{\omega}
\renewcommand\phi{\varphi}

\newcommand\chr{{\rm Chr}}
\newcommand\qed{\hfill{\vbox{\hrule\hbox{\vrule\kern3pt
                \vbox{\kern6pt}\kern3pt\vrule}\hrule}}}
\newcommand\forc{{\,\parallel\joinrel\relbar\joinrel\relbar\,}}

\begin{document}
\maketitle

\begin{abstract}
It is consistent that for every function $f:\o\to\o$ there is a graph with 
size and chromatic number $\aleph_1$ in which every $n$-chromatic subgraph 
contains at least $f(n)$ vertices ($n\geq 3$). 
This solves a \$ 250 problem of Erd\H{o}s. 
It is consistent that there is a graph $X$ with $\chr(X)=|X|=\aleph_1$ 
such that if $Y$ is a graph all whose finite subgraphs occur in 
$X$ then $\chr(Y)\leq \aleph_2$ (so the Taylor conjecture may fail). 
It is also consistent that if $X$ is a graph with chromatic number at least 
$\aleph_2$ then for every cardinal $\l$ there exists a graph $Y$ with 
$\chr(Y)\geq\l$ all whose finite subgraphs are induced subgraphs of $X$.
\end{abstract}

\section{Introduction}

In \cite{erdhaj1}  Erd\H os and Hajnal determined those finite graphs 
which appear as subgraphs in every uncountably chromatic graph: 
the bipartite graphs. 
In fact, not just that any odd circuit can be omitted, 
for every natural number $n\geq 1$ and infinite cardinal $\k$ there 
is a graph with cardinality and chromatic number $\k$ such that it 
omits all odd circuits up to length $2n+1$. 
They observed that the so-called $r$-shift graph construction has all but 
one of these properties; 
the vertex set of ${\rm Sh}_r(\kappa)$ is the set of all $r$-tuples 
from $\kappa$, with $\{x_0,x_1,\dots,x_{r-1}\}_<$ joined to 
$\{x_1,x_2,\dots,x_{r}\}_<$, then this graph omits odd circuits 
of length $3,5,\dots,2r+1$ 
and  the Erd\H os-Rado theorem asserts that the chromatic number of 
${\rm Sh}_r({\rm exp}_{r-1}(\kappa)^+)$ is at least $\kappa^+$.

The problem of determining the classes of finite graphs that occur 
in uncountably chromatic graphs seems to be much harder, and its 
investigation  was strongly 
pushed by Erd\H os and Hajnal. 

An early conjecture for example was the following. 
Every uncountably chromatic graph contains all odd circuits from some length 
onward.
This was then proved by Erd\H os, Hajnal, and Shelah \cite{erdhajshe}, 
and independently, by Thomassen \cite{thomassen}. 

In \cite{taylor2} and later in \cite{erdhajshe} the following problem was 
posed (the Taylor conjecture). 
If $\kappa$, $\lambda$ are uncountable cardinals and $X$ is a 
$\kappa$-chromatic graph, is there a $\lambda$-chromatic graph $Y$ such that every finite
subgraph of $Y$ appears as a subgraph of $X$. 
Notice that the above shift graphs give some evidence for this 
conjecture---the finite subgraphs of ${\rm Sh}_r(\kappa)$ do not depend 
on the parameter $\k$. 
The authors of \cite{erdhajshe} remarked that even the following 
much stronger conjecture seemed possible. 
If $X$ is uncountably chromatic, then for some  $r$ it contains all finite
subgraphs of ${\rm Sh}_r(\omega)$. 
This conjecture was then disproved in \cite{hajkom}. 

One easy remark, the so called Hanf number argument gives that 
there is a cardinal $\k$ with the property that if the chromatic number 
of some graph $X$ is at least $\k$ then there are arbitrarily large 
chromatic graphs with all finite subgraphs appearing in $X$. 
This argument, however, does not give any reasonable bound on $\k$. 

Another conjecture of Erd\H os and Hajnal if the maximal chromatic number of 
$n$-element subgraphs of an uncountably chromatic graph as a function of $n$ 
can converge to infinity arbitarily slowly as $n$ tends to infinity. 
It was mentioned in several problem papers, for example in \cite{erd1}, 
\cite{erd4}, \cite{erd5}, \cite{erdhaj3}, \cite{erdhajsze}. 
See also \cite{chunggraham}, \cite{jensentoft}. 
The relevance of the above examples is that the chromatic number of  the 
$n$-vertex subgraphs of (any) ${\rm Sh}_r(\kappa)$ grows roughly as 
the $r-1$ times iterated logarithm of $n$. 
Perhaps it was this fact that led Erd\H os and Hajnal to the above problem.

Erd\H os also tirelessly popularized the Taylor conjecture, he mentioned it 
e.g., in  \cite{erd2}, \cite{erd3}, \cite{erd6}, \cite{erdhaj2}, 
\cite{erdhaj3}, \cite{erdhajshe}. 
It is also mentioned in \cite{chunggraham}, the book collecting Erd\H os' 
conjectures on graphs. 
In \cite{komshe}  we gave some results when the
additional hypotheses $\vert X\vert =\kappa$, $\vert Y\vert
=\lambda$ was imposed. 
We described  countably many different classes  ${\cal K}_{n,e}$ of finite
graphs and proved that if $\lambda^{\aleph_0}=\lambda$ then every
$\lambda^+$-chromatic graph of cardinal $\lambda^+$ contains, 
for some $n$, $e$, all elements of ${\cal K}_{n,e}$ as subgraphs. 
On the other hand, it is consistent for
every regular infinite cardinal $\kappa$ that there is a
$\kappa^+$-chromatic graph on $\kappa^+$ that contains finite
subgraphs only from ${\cal K}_{n,e}$. 
We got, therefore, some models of set
theory, where the finite subraphs of graphs with 
$\vert X\vert=\chr(X) =\kappa^+$ for regular uncountable cardinals $\kappa$
were completely described. 

Notice that the class of regular cardinals on which the above result 
operated excludes $\o_1$, and in this paper we show the reason, 
by resolving the above Erd\H os-Hajnal conjecture: 
it is consistent that for every monotonically increasing  function 
$f:\o\to\o$ there is a graph with size and chromatic number 
$\aleph_1$ in which every $n$-chromatic subgraph has at least 
$f(n)$ elements ($n\geq 3$). 
The possibility of transforming the proof into a ZFC argument will be checked 
in the forthcoming \cite{Sh:e}, Chapter 9. 
An application of the  method presented here gives the consistent existence 
of a graph $X$ with $\chr(X)=|X|=\aleph_1$ such that if $Y$ is a graph 
(of any size) 
all whose subgraphs are subgraphs of $X$ then $\chr(Y)\leq\aleph_2$. 
This gives a consistent negative answer to the Taylor conjecture. 
As for the positive direction we prove that it is consistent that 
if $X$ is a graph with chromatic number at least $\aleph_2$ then there are 
arbitrarily large chromatic graphs all whose finite subgraphs 
being induced subgraphs of $X$.

Theorems 1 and 2 were proved by S.~Shelah and then P.~Komj\'ath 
proved Theorems 3 and 4. 

\medskip
\noindent{\bf Notation.} We use the standard axiomatic set theory notation. 
If $A$ is a set of ordinals, $\b$ is an ordinal, then 
$\b<A$,  means that $\b<\a$ holds for every $\a\in A$. 
Similarly for $\b\leq A$, $A<\b$, etc. 
If $f$ is a function, $A$ a set, then we let $f[A]=\{f(x):x\in A\}$. 
If $S$ is a set, $\kappa$ is a cardinal,  
$[S]^\kappa=\{X\subseteq S:|X|=\kappa\}$, 
$[S]^{<\kappa}=\{X\subseteq S:|X|<\kappa\}$.  
A {\sl graph} is an ordered pair $(V,X)$ where $V$ is some set 
(the set of {\em vertices}) 
and $X\subseteq [V]^2$ (the set of {\sl edges}). 
In some cases we identify the graph and $X$. 
The {\em chromatic number} of a graph $(V,X)$ is the least cardinal 
$\mu$ such that there exists a function $f:V\to\mu$ with $f(x)\neq f(y)$ 
for $\{x,y\}\in X$ (a good coloring). 
A {\em path} of length $n$ is a sequence $\{v_0,\dots,v_n\}$ 
of distinct vertices such that $\{v_i,v_{i+1}\}\in X$ holds for 
$i<n$.
A  {\em cycle} of length $n$ is a sequence $\{v_1,\dots,v_n\}$ 
of vertices such that for $1\leq i<n$ we have 
$\{v_i,v_{i+1}\}\in X$ and $\{v_n,v_{1}\}\in X$ also holds. 
If the vertices are distinct then we call it a {\em circuit}.

\section{A large chromatic graph with small chromatic finite subgraphs} 
\medskip
Let $f:\o\to\o$ be a strictly increasing function. 
Fix a sequence $\{C_\a:\a<\o_1,\mbox{limit}\}$ such that $C_\a$ 
is an $\o$-sequence converging to $\a$, and the whole sequence is a club 
guessing sequence, that is, 
if $C\subseteq\o_1$ is a closed unbounded set, then 
$C_\a\subseteq C$ holds for some $\a$. 
Notice that the existence of a guessing sequence is an easy consequence 
of the diamond principle. 

We are going to define the notion of forcing $(Q^f,\leq)$.
Every condition $p\in Q^f$ will be of the form 
$p=(s,u,g,m,h_i,c)$ where $s\in [\o_1]^{<\o}$, 
$u\subseteq s$, consisting of limit ordinals, $g$ is a graph on $u$, 
for every $\a\in u$ we have an $m(\a)<\o$, and then the ordinals 
$h_0(\a)<\cdots<h_{m(\a)-1}(\a)<\a$ which are $C_\a$-separated, that is, 
$\min(C_\a)<h_0(\a)$ and between $h_{i}(\a)$ and $h_{i+1}(\a)$ there is an 
element of $C_\a$ ($h_i(\a)$ is undefined if $i\geq m(\a)$ or 
$\a\notin u$). 
$c:g\to\o$ satisfies that if $\{x,\a\}\in g$, $\{y,\a\}\in g$ and 
$x,y<\a$ then $c(x,\a)\neq c(y,\a)$. 

Given $p\in Q^f$ we define 
\[ 
         Y^p_r = \{e\in g:c(e)\geq r\} 
\]
for any natural number $r$.         

We add the following stipulations. 

\begin{enumerate} 
\item[(1.)] If $\a\in u$ is incident to some $e\in g$ with $c(e)\geq r$ 
          then $m(\a)\geq r$. 
\item[(2.)] If $x$, $y$ are connected in $Y^p_r$ in $i$ steps then 
          $h_{j-i}(x)<h_{j}(y)<h_{j+i}(x)$ holds for $i\leq j\leq f(r)-i$. 
\item[(3.)] $Y^p_r$ does not contain odd circuits of length $\leq f(r)$. 
\end{enumerate}

We notice that although a condition 
$p=(s,u,g,m,h_i,c)$ is infinite (perhaps   
$p=(s,u,g,m,(h_i)_{i<\o},c)$ or $p=(s,u,g,m,h_i,c)_{i<\o}$ would be 
better notation) 
as all but finitely many of the partial functions $\{h_i:i<\o\}$ 
are the empty function, every condition is really a finite object. 

The partial order on $Q^f$ is defined the natural way. 
$p'=(s',u',g',m',h'_i,c')$ extends $p=(s,u,g,m,h_i,c)$ iff the following 
hold. 
$s'\supseteq s$, $u=u'\cap s$, $g=g'\cap [s]^2$. 
$m'(\a)\geq m(\a)$ 
holds for $\a\in u$, and 
$h'_i(\a)=h_i(\a)$ for $i<m(\a)$. 
Finally, $c'(x,\a)=c(x,\a)$ holds for $\{x,\a\}\in g$. 

We call two conditions $p=(s,u,g,m,h_i,c)$ and $p'=(s',u',g',m',h'_i,c')$ 
{\em isomorphic} iff $|s|=|s'|$ and the unique order preserving mapping 
$\pi:s\to s'$ satisfies $u'=\pi[u]$, $g'=\pi[g]$, 
$m(\a)=m'(\pi(\a))$ for $\a\in u$, $h'_i(\pi(\a))=h_i(\a)$ and 
$c'(\pi(x),\pi(\a))=c(x,\a)$ hold whenever the right hand sides are defined. 
Notice that, as our conditions are finite structures, we have only 
countably many isomorphism types.

{}From a generic $G\subseteq Q^f$  we define the 
following graphs on $\o_1$: 
\[ 
   X = \bigcup \left\{g:(s,u,g,m,h_i,c)\in G \right\}, 
\]
\[ 
   X_n =  \left\{ \{x,y\}: \{x,y\}\in g , c(x,y)=n,(s,u,g,m,h_i,c)\in G
          \right\}, 
\]

\[ 
   Y_r = X_r\cup X_{r+1}\cdots 
\]
and notice that $Y_r=\cup\{Y^p_r:p\in G\}$.

\medskip
\noindent{\bf Lemma 1.} 
         {\em For $\a<\o_1$ the set 
         $\{(s,u,g,m,h_i,c)\in Q^f:\a\in s\}$ is dense.}

\medskip
\noindent{\bf Proof.} 
Straightforward.
\qed

\medskip
\noindent{\bf Lemma 2.} 
         {\em If $\a<\o_1$ is limit, $(s,u,g,m,h_i,c)\in Q^f$, $\a\notin s$, 
         then there is an extension 
         $(s',u',g',m',h'_i,c')\leq (s,u,g,m,h_i,c)$ with $\a\in u'$.}

\medskip
\noindent{\bf Proof.} 
Straightforward.
\qed

\medskip
\noindent{\bf Lemma 3.} 
         {\em If $(s,u,g,m,h_i,c)\in Q^f$, 
         $\a\in u$, $n<\o$ then there is an extension 
         $(s',u',g',m',h'_i,c')\leq (s,u,g,m,h_i,c)$ with $m(\a)\geq n$.}

\medskip
\noindent{\bf Proof.} 
It suffices to show that $m(\a)$ can be incremented by one. 
Given $h_{m(\a)-1}(\a)<\a$ if we choose $h'_{m(\a)}(\a)<\a$ large enough 
the condition on $C_\a$ will surely be satisfied. 
\qed 

\medskip
\noindent{\bf Lemma 4.} 
         {\em $(Q^f,\leq)$ is ccc.}

\medskip
\noindent{\bf Proof.} 
Modulo standard arguments we have to show that $p=(s,u,g,m,h_i,c)$ 
and $p'=(s',u',g',m',h'_i,c')$ are compatible, assuming that they are 
isomorphic, and $s\cap s'$ is an initial segment of both $s$ and $s'$. 
Let $\pi:s\to s'$ be the order preserving structure isomorphism. 
We let $p''=(s'',u'',g'',m'',h''_i,c'')$ where we take unions in all 
coordinates. 

In order to show that $p''$ is a condition we have to check 
properties (1.--3.). 

(1.) is obvious. 

For (3.) assume that $C$ is an odd circuit of length $\leq f(r)$ 
in $Y^{p''}_r$. 
If we replace every $e\in C$ that contains at least one vertex from $s-s'$ 
with $\pi(e)$ then we get an odd cycle $C'$ in $Y^{p'}_r$. 
$C'$ splits into circuits, at least one of them odd, so we get a 
contradiction. 

For (2.) notice that it holds if $x$, $y$ are joined via a path 
going entirely in $Y^{p}_{r}$ or $Y^{p'}_{r}$. 
It suffices, therefore, to show that if some path $P$ between $x$ and $y$ 
of length $i$ is split by an inner point $z$ into the paths 
$P_0$ and $P_1$ between $x$ and $z$, and $z$ and $y$, respectively, 
and of the respective lengths $i_0$ and $i_1$ (so $i_0+i_1=i$) 
and the statement holds for $P_0$ and $P_1$ then it holds for $P$, as well. 
Indeed, for $i\leq j \leq f(r)-i$ we have 
\[
h_{j-(i_0+i_1)}(x)<h_{j-i_1}(z)<h_{j}(y)<h_{j+i_1}(z)<h_{j+(i_0+i_1)}(x).
\]
\qed

\medskip
\noindent{\bf Lemma 5.} 
         {\em $\chr(X)=\o_1$.}

\medskip
\noindent{\bf Proof.} 
Assume that some $\overline{p}$ forces that $\underline{F}$ is a good 
coloring of $X$ with the elements of $\o$.  
Select an increasing, continuous sequence of countable elementary submodels 
$ \overline{p}, \underline{F}, Q^f, 
\forc\in N_0\prec N_1\prec\cdots N_\a \prec H(\lambda)$ with some large 
enough regular cardinal $\lambda$, for $\a<\o_1$, such that 
$\g_\a=N_\a\cap\o_1$ is an ordinal. 
The set $C=\{\g_\a:\a<\o_1\}$ will be a closed, unbounded set, so by the 
guessing property there is some $\d=\g_\a$ such that $C_\d\subseteq C$. 
Notice that all points of $\overline{p}$ are smaller than $\d$. 

Extend $\overline{p}$ to a $p'$ using Lemma 2., adding $\d$ to the 
$u$-part, 
then let $p^*$ be a condition extending $p'$ such that 
$p^*\forc \underline{F}(\d)=i$ holds for some $i<\o$.

Let $n$ be some natural number that 
$n\notin \{c^*(x,\d):\{x,\d\}\in g^*, x<\d\}$ where 
$p^*=(s^*,u^*,g^*,m^*,h^*_i,c^*)$. 
Extend $p^*$ using Lemma 3., to some condition $p$  with 
$p=(s,u,g,m,h_i,c)$ such that $m=m(\d)\geq f(n)$ holds. 

In $p$ we have the values 
\[ 
      h_{0}(\d)<h_{1}(\d)<\cdots<h_{m-1}(\d)<\d
\]
and by our requirements on conditions there are elements 
$\d_0,\dots,\d_{m-1}$ of $C_\d$ such that 

\[ 
      \d_0<h_{0}(\d)<\d_1<h_{1}(\d)<\cdots<\d_{m-1}<h_{m-1}(\d)<\d
\]
holds. 

The values $\d_0,\dots,\d_{m-1},\d$ break $s$ into disjoint parts:       
$s=s_0\cup\cdots\cup s_m\cup s_{m+1}$ 
with 
\[
s_0<\d_0\leq s_1<\d_1\leq s_2<\cdots <\d_{m-1}\leq s_m <\d\leq s_{m+1}
\]
(some of them may be empty).

\medskip
\noindent{\bf Sublemma.} 
         {\em There is a condition $p'$ on 
         some $s'=s_0\cup s'_1\cup\cdots\cup s'_m\cup s'_{m+1}$ 
         isomorphic to $p$ with $|s'_i|=|s_i|$,   
         $\d'=\min(s'_{m+1})$, $p'\forc \underline{F}(\d')=i$ and 
         \[
         s_0<s'_1<\d_0\leq s_1<s'_2<\d_1\leq s_2<\cdots< s'_m
         <\d_{m-1}\leq s_m <s'_{m+1} <\d\leq s_{m+1}
         \]
         holds.}

\medskip
\noindent{\bf Proof.} 
Let $\th$ be the isomorphism type of $p$. 
For the ordered finite sets 
$\overline{x}_0,\overline{x}_1,\dots,\overline{x}_{m+1}$
let $\psi(\overline{x}_0,\overline{x}_1,\dots,\overline{x}_{m+1})$ 
denote the statement that 
$\overline{x}_0<\overline{x}_1<\cdots<\overline{x}_{m+1}$, 
$|\overline{x}_i|=|s_i|$ and for the 
(unique) condition $p$ on 
$\overline{x}_0\cup\overline{x}_1\cup\cdots\cup\overline{x}_{m+1}$ of type 
$\th$ 
$p\forc \underline{F}(\d)=i$ where $\d=\min(\overline{x}_{m+1})$. 

Set 
\[
\varphi_{m+1}(\overline{x}_0,\overline{x}_1,\dots,\overline{x}_{m+1})=
\psi(\overline{x}_0,\overline{x}_1,\dots,\overline{x}_{m+1})
\] 
and for 
$0\leq i\leq m$ 
\[
\varphi_i(\overline{x}_0,\overline{x}_1,\dots,\overline{x}_{i})=
\exists^*\overline{x}_{i+1}\varphi_{i+1}(\overline{x}_0,\overline{x}_1,
\dots,\overline{x}_{i+1})
\]
where the quantifier $\exists^*$ denotes ``there exist unboundedly many'' 
which is expressible in the first order language of $(\o_1,<)$.  

\medskip
\noindent{\bf Claim.} 
         {\em For $0\leq i \leq m+1$ the sentence 
         $\varphi_i(s_0,s_1,\dots,s_{i})$ holds.}

\medskip
\noindent{\bf Proof.} 
We prove the statement by reverse induction on 
$0\leq i \leq m+1$. 
We certainly have $\varphi_{m+1}(s_0,s_1,\dots,s_{m+1})$. 
If for some $0\leq i \leq m$ we had that 
$\varphi_{i+1}(s_0,s_1,\dots,s_{i+1})$ holds yet 
\[
\varphi_i(s_0,s_1,\dots,s_{i})=
\exists^*\overline{x}_{i+1}\varphi_{i+1}(s_0,s_1,\dots,s_i,
\overline{x}_{i+1})
\] 
fails, then there was a bound, computable from $s_0,s_1,\dots,s_i$ for 
the minima of those sets $\overline{x}_{i+1}$ for which 
$\varphi_{i+1}(s_0,s_1,\dots,s_i,\overline{x}_{i+1})$ holds. 
Then this bound was smaller than $\d_i$ ($\d$ for $i=m$) as there 
is an elementary submodel containing the ordinals $<\d_i$ (or $<\d$) 
but this contradicts the fact that $\varphi_{i+1}(s_0,s_1,\dots,s_{i+1})$ holds 
and $\d_i\leq s_i$. 
\qed 

\medskip
Using the Claim we can inductively select the sets $s'_1,\dots,s'_{m+1}$ such 
that for every $1\leq i \leq m+1$ we have $\varphi_i(s_0,s'_1,\dots,s'_{i})$ 
and $s_i<s'_{i+1}\leq \d_i$, as required. 
\qed

\medskip
Using the Sublemma we create the following one-edge amalgamation $p''$ of 
$p$ and $p'$. 
$p''=(s'',u'',g'',m'',h''_i,c'')$ where 
$s''=s\cup s'$,  
$u''=u\cup u'$,  
$g''=g\cup g'\cup\left\{\{\d,\d'\}\right\}$, 
$m''=m\cup m'$, 
$h''_i=h_i\cup h'_i$, 
$c''=c\cup c'$. 

We have to show that $p''$ is indeed a condition, that is, we have to check 
if the properties (1.-3.) hold. 

(1.) is obvious. 

For (2.) we argue as in the proof of Lemma 4; every path in question 
is the union of paths for which this condition holds, plus 
possibly the path $\{\d,\d'\}$ but (2.) also holds for this. 

Assume finally, that $C$ is a circuit of length $2t+1\leq f(r)$ in 
$Y^{p''}_r$. 
Unless $C$ contains $\{\d,\d'\}$, we can argue as in Lemma 4. 
So we are left with the case that $C$ contains $\{\d,\d'\}$ and 
$r\leq n$. 
That is, $\d \in s$ and $\d' \in s'$ are joined in 
$Y^{p}_r\cup Y^{p'}_r$ in $2t$ steps, and this is only possible 
if the connecting path has vertices  in $s\cap s'$. 
So we get that $\d$ can be connected in $Y^{p}_r$ with some point in 
$s\cap s'$ in $\leq t$ steps. 
But this is impossible: if $x\in s\cap s'$ is such a point then 
$ h_0(\d)<h_t(x)<h_{2t}(\d)$ by condition (2.) and also $x<h_0(\d)$, a 
contradiction.  

As $p''$ forces that $F(\d)=F(\d')=i$ yet they are joined in $X$, 
we are finished. 
\qed

\medskip
\noindent{\bf Theorem 1.} 
         {\em The forcing $Q^f$ adds an uncountably chromatic graph $X$ 
         on $\o_1$ such that every subgraph on at most $f(r)$ vertices 
         is at most $2^{r+1}$-chromatic.}

\medskip
\noindent{\bf Proof.} 
As every $X_n$ is a circuitfree graph, it can be colored with two colors. 
Consider now a subgraph of $X$ induced by a set $S$ of at most $f(r)$ 
vertices. 
On $S$ all the graphs $X_0,\dots,X_{r-1}$ are bipartite, and so is 
$Y_r=X_r\cup\cdots$ (as it has no odd circuits of length $\leq f(r)$). 
So their union, $X$ restricted to $S$, can be colored by at most $2^{r+1}$ 
colors. 
\qed 

\medskip
\noindent{\bf Theorem 2.} 
         {\em It is consistent with CH that for every function $f:\o\to\o$
         there is an uncountably 
         chromatic  graph $X$ on $\o_1$ such that  
         every sugraph of $X$ on $f(r)$ vertices 
         is at most $r$-chromatic ($r\geq 2$).}

\medskip
\noindent{\bf Proof.} 
Assume that $\diamondsuit$ holds in the ground model. 
Then we have CH and there is a club guessing sequence 
$\{C_\a:\a<\o_1,\mbox{limit}\}$ as required for 
Theorem 1. 
We force with a finite support iteration $P=\{P_\a,Q_\a:\a<\o_1\}$. 
In step $\a<\o_1$ we add $Q_\a=Q^{f_\a}$ for some increasing function 
$f_\a:\o\to\o$. 
As this will be a ccc forcing that preserves CH it is possible by 
bookkeeping to make sure that every suitable $f:\o\to\o$ occurs as 
some $f_\a$. 
Also, as $P_\a$, the iteration up to $\a$ is ccc, every closed, unbounded 
set $C$ in $V^{P_\a}$ contains a ground model closed, unbounded set $D$, and 
as there is some element $C_\a$ of the club guessing system that 
$C_\a\subseteq D$ we have $C_\a\subseteq C$, that is the club 
guessing system retains its property in $V^{P_\a}$.

Call a condition $p\in P$ {\em determined} 
if for every 
$\a<\o_1$ the condition $p|\a$ completely determines $p(\a)$, 
that is, for every $\a$  coordinate $p(\a)$ is not just a name 
for a finite structure but it is actually a finite structure. 

\medskip
\noindent{\bf Lemma 6.} 
         {\em The determined conditions form a dense set in $P$.}

\medskip
\noindent{\bf Proof.} 
We prove by induction on $\a<\o_1$ that the determined conditions form 
a dense subset of $P_\a$. 
This is obvious if $\a$ is limit, as we are considering finite supports. 
Assume that we have the statement for some $\a<\o_1$ and we try to 
handle the case of $\a+1$. 
Let $(p,q)\in P_{\a+1}=P_\a\ast Q_\a$ be arbitrary. 
Extend $p$ to some $p'$ that completely determines $q$, that is, there is 
a finite structure $h$ that $p'\forc q=h$. 
Then extend $p'$ to a determined $p^*\in P_\a$. 
Now $(p^*,h)$ is a determined extension of $(p,q)$. 
\qed

\medskip
\noindent{\bf Lemma 7.} 
         {\em For every $\a<\o_1$, $\chr(X_\a)=\o_1$ holds in $V^P$.}

\medskip
\noindent{\bf Proof.} 
By moving to $V^{P_\a}$ we can assume that $\a=0$. 
We imitate the proof of Lemma 5. 
By Lemma 6 we can work with determined conditions. 
We consider every such condition as a finite structure on some finite 
subset $s$ of $\omega_1$, here $s$ contains all points of all graphs 
$p(\b)$ where $\b$ is an arbitrary element of the support of $p$, 
and we also add the elements of the support to $s$. 
Assume that some $\overline{p}\in P$ forces that $\underline{F}$ is a good 
coloring of $X_0$ with the elements of $\o$. 
With an argument like in Lemma 5 we get some natural number $n$, ordinals 
$\d'<\d<\o_1$, and also $\a_1,\dots,\a_m,\b_1,\dots,\b_t, \b'_1,\dots,\b'_t$ 
and two isomorphic determined conditions $p$ and $p'$ with the respective 
supports $\{0,\a_1,\dots,\a_m,\b_1,\dots,\b_t\}$ and 
$\{0,\a_1,\dots,\a_m,\b'_1,\dots,\b'_t\}$ that $p\forc\underline{F}(\d)=n$, 
$p'\forc\underline{F}(\d')=n$ hold, $p(0)$ and $p'(0)$ are isomorphic 
conditions that behave like $p$ and $p'$ in Lemma 5, for every $\a_i$ 
the structures $p(\a_i)$ and $p'(\a_i)$ are isomorphic with the 
common part preceding the tails, that is, it is possible to make 
a non-edge amalgamation. 
We can, therefore, take the union of $p$ and $p'$, and add the edge 
$\{\d,\d'\}$ with color $n$ in coordinate 0.
\qed

\medskip
The above Lemma concludes the proof of the Theorem.
\qed

\section{The Taylor conjecture}

\medskip
\noindent{\bf Theorem 3.} 
         {\em It is consistent that there is a graph $X$ with 
         $\chr(X)=|X|=\aleph_1$ such that if $Y$ is a graph with all finite 
         subgraphs occurring in $X$ then $\chr(Y)\leq \aleph_2$, that is, 
         the Taylor conjecture fails.}

\medskip
\noindent{\bf Proof.} 
Let $V$ be a model of GCH and $\diamondsuit$. 
Let $P$ be the notion of forcing that adds a Cohen real. 
It is well known that $P$ adds an undominated real, 
that is, a function $f:\o\to\o$ such that for no $g:\o\to\o$ 
in $V$ does  $f(n)\leq g(n)$ hold for every $n<\o$. 
Let $V'$ be the forced model. 
Notice that $V'$ still has GCH and club guessing (by the argument in the 
proof of Theorem 2.).  
Now force over $V'$ with the partial order $Q^f$, and get a graph $X$ 
with $\chr(X)=|X|=\aleph_1$ such that every $n$-chromatic subgraph of 
$X$ has at least $f(n)$ elements ($n\geq 3$). 
This $X$ will be our graph. 
To show the property stated, assume that $Y$ is a graph in $V^{P,Q^f}$ 
whose every finite subgraph is a subgraph of $X$. 
We assume that the vertex set of $Y$ is some cardinal $\lambda$.
We notice that every $n$-chromatic subgraph of $Y$ has at least $f(n)$ 
elements. 

\medskip
\noindent{\bf Lemma 8.} 
         {\em If $Z\subseteq Y$ is a subgraph with $Z\in V$ then 
         $Z$ is finitely chromatic.} 

\medskip
\noindent{\bf Proof.} 
Otherwise for every $n<\o$ we can let $g(n)$ be the minimal size of 
an $n$-chromatic subgraph of $Z$. 
Now notice that $g\in V$ and also by the absoluteness of the set of 
finite subsets of $\lambda$ and the absoluteness of the cromatic number 
of a finite graph, $g$ denotes the same thing in $V$ and $V^{P\ast Q^f}$.
This implies that $g(n)\geq f(n)$ holds for every $n$,  but that obviously 
contradicts the fact that $f$ cannot be dominated by the the ground model 
$\o\to\o$ functions. 
\qed

\medskip
We finally need the following Lemma. 

\medskip
\noindent{\bf Lemma 9.} 
         {\em If $R$ is a notion of forcing over some model $V$, 
         $Y$ is a graph in the extended model on some ordinal $\l$ then $Y$ 
         is the union of at most $|R|$ subgraphs which are elements of $V$.}

\medskip
\noindent{\bf Proof.} 
Let $\tau$ be a name  for (the edge set of) $Y$. 
Set 
\[
Z_p=\{e\in [\lambda]^2:p\forc e\in\tau\}
\] 
for $p\in R$, then 
$Y=\bigcup\{Z_p:p\in G\}$ where $G\subseteq R$ is a generic set. 
\qed 

\medskip
To finish the proof of the Theorem we remark that by Lemma 9 $Y$ 
decomposes into the union of $|P\ast Q^f|=\aleph_1$ graphs 
each being in $V$, therefore finitely chromatic, so we get 
$\chr(Y)\leq 2^{\aleph_1}=\aleph_2$.
\qed

\medskip
The followig argument gives that the Hanf number mentioned in the 
Introduction can be as small as $\aleph_2$.

\medskip
\noindent{\bf Theorem 4.} 
         {\em It is consistent that if $X$ is a graph with 
         $\chr(X)\geq\aleph_2$ then for every cardinal $\l$ there exists 
         a graph $Y$ with $\chr(Y)\geq\l$ all whose finite subgraphs 
         are induced subgraphs of $X$.}

\medskip
\noindent{\bf Proof.} 
Let $V$ be a model of GCH. 
Choose the regular cardinal $\k$ so large that the following holds. 
If $X$ is  a graph with $\chr(X)\geq\k$ and $\l$ is a cardinal 
then there is a graph with $\chr(Y)\geq\l$ all whose finite 
subgraphs occur as subgraphs of $X$. 
Clearly, such a $\k$ exists. 

Let $P={\rm Col}(\o,\k)$ be the collapse of $\k$ to $\aleph_0$, 
that is, the elements of $P$ are those functions of the form 
$p:n\to\k$ for some $n<\o$ with $p\leq q$ iff $p$ extends $q$ as 
a function. 
Our claim is that if $G\subseteq P$ is generic then $V[G]$ models 
the statement of the Theorem. 
Notice that $|\k|=\aleph_0$ holds there and calculation shows that 
GCH still holds in $V[G]$. 

Assume that $X$ is a graph in $V[G]$ with chromatic number at 
least $\aleph_2$ (that is, $\aleph_2^{V[G]}$). 
By Lemma 9., $X$ is the union of $|G|=\aleph_0$ ground model graphs. 
As $\aleph_1^{\aleph_0}=\aleph_1$, one of them, say $Y$ must have 
chromatic number at least $\aleph_2$. 
In $V$, $Y$ has chromatic number at least $\aleph_2^{V[G]}=\k^{++}$. 
Assume that we are given some $\l>\k$.  
By the choice of $\k$, there is a graph $Z$ with $\chr(Z)\geq\l$ such that 
every finite induced subgraph of $Z$ is an induced subgraph of $Y$.

\medskip
\noindent{\bf Lemma 10.} 
         {\em $\chr(Z)\geq\l$ holds in $V[G]$.}

\medskip
\noindent{\bf Proof.} 
Otherwise let $\underline{F}$ be a name for a coloring with 
the ordinals less than $\tau<\l$. 
Then the coloring $x\mapsto (p,\xi)$ is a good coloring of the 
vertices of $Z$ with $\k+\tau<\l$ colors, where  $p\in P$ 
is some element of $P$ with $p\forc \underline{F}(x)=\xi$. 
\qed 

\medskip
We are almost finished, the only problem is that the finite induced 
subgraphs of $Z$ are not induced subgraphs of $X$, they only are 
(edge-)subgraphs of induced subgraphs of $X$. 
The following Lemma is what we need. 

\medskip
\noindent{\bf Lemma 11.} 
         {\em There is a graph $Z'$ on the vertex set of $Z$ with 
         $Z'\supseteq Z$ and such that every induced subgraph of $Z'$ is an 
         induced subgraph of $X$.}

\medskip
\noindent{\bf Proof.} 
Let $S$ be the vertex set of $Z$. 
For every finite subset $s$ of $S$ there are some graphs on $s$,  
which on the one hand are isomorphic to induced subgraphs 
of $X$, on the other hand they are supergraphs of the graph 
$Z$ restricted to $s$. 
Call these graphs {\em appropriate} for $s$. 
Notice that there are finitely many appropriate graphs for every given $s$, 
and if $T$ is  an appropriate graph for $s$ and $s'$ is a subset of $s$ then 
$T$ restricted to $s'$ is a graph appropriate for $s'$. 
We can therefore apply the Rado selection principle 
(or the compactness theorem of model theory) and get a graph 
$Z'$ on $S$ every induced subgraph of which is appropriate, so 
$Z'$ is a required. 
\qed 

As $\chr(Z')\geq\chr(Z)\geq\l$ holds we are done. 
\qed

\vskip1cm
\hbox{
\vtop{\hbox{P\'eter Komj\'{a}th}
      \hbox{Department of Computer Science}
      \hbox{E\"otv\"os University}
      \hbox{Budapest, P.O.Box 120 } 
      \hbox{1518, Hungary}
      \hbox{e-mail:{\tt\ kope@cs.elte.hu}}}
      \qquad\qquad\quad
\vtop{\hbox{Saharon Shelah}
      \hbox{Institute of Mathematics,} 
      \hbox{Hebrew University,} 
      \hbox{Givat Ram, 91904,} 
      \hbox{Jerusalem, Israel}
      \hbox{e-mail:{\tt\ shelah@math.huji.ac.il}}}}

\end{document}